\documentclass[fleqn]{llncs}
\usepackage{am}
\usepackage{version}

\title{Arithmetical Meadows}
\author{J.A. Bergstra \and C.A. Middelburg}
\institute{Informatics Institute, Faculty of Science,
           University of Amsterdam, \\
           Science Park~904, 1098~XH Amsterdam, the Netherlands \\
           \email{J.A.Bergstra@uva.nl,C.A.Middelburg@uva.nl}}

\begin{document}

\maketitle

\begin{abstract}
An inversive meadow is a commutative ring with identity equipped with a
total multiplicative inverse operation satisfying $0\minv = 0$.
Previously, inversive meadows were shortly called meadows.
A divisive meadow is an inversive meadows with the multiplicative
inverse operation replaced by a division operation.
In the spirit of Peacock's arithmetical algebra, we introduce variants
of inversive and divisive meadows without an additive identity element
and an additive inverse operation.
We give equational axiomatizations of several classes of such variants
of inversive and divisive meadows as well as of several instances of
them.
\begin{keywords}
arithmetical meadow, equational specification, initial algebra
specification.
\end{keywords}%
\begin{classcode}
12E12, 12E30, 12L05.
\end{classcode}
\end{abstract}

\section{Introduction}
\label{sect-introduction}

The primary mathematical structure for measurement and computation is
unquestionably a field.
In~\cite{BT07a}, meadows are proposed as alternatives for fields with a
purely equational specification.
A meadow is a commutative ring with identity equipped with a total
multiplicative inverse operation satisfying two equations which imply
that the multiplicative inverse of zero is zero.
Thus, meadows are total algebras.
As usual in field theory, the convention to consider $p \mdiv q$ as an
abbreviation for $p \mmul q\minv$ was used in subsequent work on
meadows (see e.g.~\cite{BHT09a,BP08a}). \linebreak[2]
This convention is no longer satisfactory if partial variants of meadows
are considered too, as is demonstrated in~\cite{BM09g}.
In that paper, we rename meadows into inversive meadows and introduce
divisive meadows.
A divisive meadow is an inversive meadow with the multiplicative inverse
operation replaced by the division operation suggested by the
above-mentioned abbreviation convention.
Henceforth, we will use the name meadow whenever the distinction between
inverse meadows and divisive meadows is not important.

Peacock introduced in~\cite{Pea1830a} arithmetical algebra as algebra of
numbers where an additive identity element and an additive inverse
operation are not involved.
That is, arithmetical algebra is algebra of positive numbers instead of
algebra of numbers in general (see also~\cite{Kle98a}).
In the spirit of Peacock, we use the name \emph{arithmetical meadow} for
a meadow without an additive identity element and an additive inverse
operation.
We use the name \emph{arithmetical meadow with zero} for a meadow
without an additive inverse operation, but with an additive identity
element.
Arithmetical meadows related to the field of rational numbers are
reminiscent of Peacock's arithmetical algebra.

In this paper, we pursue the following objectives:
\begin{enumerate}
\item
to complement the signatures of inversive and divisive meadows with
arithmetical versions;
\item
to provide equational axiomatizations of several classes of arithmetical
meadows and instances of them related to the field of rational numbers;
\item
to state a number of outstanding questions concerning arithmetical
meadows.
\end{enumerate}

This paper is organized as follows.
First, we go into the background of the work presented in this paper
with the intention to clarify and motivate this work
(Section~\ref{sect-background}).
Next, we introduce the classes of inversive and divisive arithmetical
meadows and the classes of inversive and divisive arithmetical meadows
with zero (Section~\ref{sect-arith-meadows}).
After that, we have an interlude on inversive and divisive meadows
(Section~\ref{sect-meadows}).
Then, we introduce instances of the classes of arithmetical meadows and
arithmetical meadows with zero related to the field of rational numbers
(Section~\ref{sect-arith-meadows-rat}).
Following this, we state some outstanding questions about arithmetical
meadows (Section~\ref{sect-questions}).
After that, we shortly discuss partial variants of the instances of the
inversive and divisive arithmetical meadows with zero introduced before
(Section~\ref{sect-partial-arith-meadows}) and an arithmetical version
of a well-known mathematical structure closely related to inversive
meadows (Section~\ref{sect-relation-arith-rings}).
Finally, we make some concluding remarks
(Section~\ref{sect-conclusions}).

\section {Background on the Theory of Meadows}
\label{sect-background}

In this section, we go into the background of the work presented in this
paper with the intention to clarify and motivate this work.

The theory of meadows, see e.g.~\cite{BHT09a,BP08a,BM09g}, constitutes a
hybrid between the theory of abstract data type and the theory of rings
and fields, more specifically the theory of von Neumann regular
rings~\cite{McC64a,Goo79a} (all fields are von Neumann regular rings).

It is easy to see that each meadow can be reduced to a commutative von
Neumann regular ring with a multiplicative identity.
Moreover, we know from~\cite{BHT09a} that each commutative von Neumann
regular ring with a multiplicative identity can be expanded to a meadow,
and that this expansion is unique.
It is easy to show that, if $\morph{\phi}{X}{Y}$ is an epimorphism
between commutative rings with a multiplicative identity and $X$ is a
commutative von Neumann regular ring with a multiplicative identity,
than:
 (i)~$Y$ is a commutative von Neumann regular ring with a multiplicative
identity;
(ii)~$\phi$ is also an epimorphism between meadows for the meadows $X'$
and $Y'$ found by means of the unique expansions for $X$ and $Y$,
respectively.

However, there is a difference between commutative von Neumann regular
rings with a multiplicative identity and meadows: the class of all
meadows is a variety and the class of all commutative von Neumann
regular rings with a multiplicative identity is not.
In particular, the class of commutative von Neumann regular rings with a
multiplicative identity is not closed under taking subalgebras (a
property shared by all varieties).
Let $\Rat$ be the ring of rational numbers, and let $\Int$ be its
subalgebra of integers.
Then $\Rat$ is a field and for that reason a commutative von Neumann
regular ring with a multiplicative identity, but its subalgebra $\Int$
is not a commutative von Neumann regular ring with a multiplicative
identity.

In spite of the fact that meadows and commutative von Neumann regular
rings with a multiplicative identity are so close that no new
mathematics can be expected, there is a difference which matters very
much from the perspective of abstract data type specification.
$\Rat$, the ring of rational numbers, is not a minimal algebra, whereas
$\Ratzi$, the inversive meadow of rational numbers is a minimal algebra.
As such, $\Ratzi$ is amenable to initial algebra specification.
The first initial algebra specification of $\Ratzi$ is given
in~\cite{BT07a} and an improvement due to Hirshfeld is given
in~\cite{BM09g}.
When looking for an initial algebra specification of $\Rat$, adding a
total multiplicative inverse operation satisfying $0\minv = 0$ as an
auxiliary function is the most reasonable solution, assuming that a
proper constructor as an auxiliary function is acceptable.

We see a theory of meadows having two roles:
 (i)~a starting-point of a theory of mathematical data types;
(ii)~an intermediate between algebra and logic.

On investigation of mathematical data types, known countable
mathematical structures will be equipped with operations to obtain
minimal algebras and specification properties of those minimal algebras
will be investigated.%
If countable minimal algebras can be classified as either computable,
semi-computable or co-semi-computable, known specification techniques
may be applied (see~\cite{BT95a} for a survey of this matter).
Otherwise data type specification in its original forms cannot be
applied.
Further, one may study $\omega$-completeness of specifications and term
rewriting system related properties.%

It is not a common viewpoint in algebra or in mathematics at large that
giving a name to an operation, which is included in a signature, is a
very significant step by itself.
However, the answer to the notorious question ``what is $1 \mdiv 0$'' is
very sensitive to exactly this matter.
Von Neumann regular rings provide a classical mathematical perspective
on rings and fields, where multiplicative inverse (or division) is only
used when its use is clearly justified and puzzling uses are rejected as
a matter of principle.
Meadows provide a more logical perspective to von Neumann regular rings
in which justified and unjustified use of multiplicative inverse cannot
be easily distinguished beforehand.

Arithmetical meadows, in particular the ones related to the field of
rational numbers, provide additional insight in what is yielded by the
presence of an operator for multiplicative inverse (or division) in a
signature.

\section{Equational Specifications of Arithmetical Meadows}
\label{sect-arith-meadows}

This section concerns the equational specification of several classes of
arithmetical meadows.

The signature of commutative rings with a multiplicative identity
consists of the following constants and operators:
\begin{itemize}
\item
the \emph{additive identity} constant $0$;
\item
the \emph{multiplicative identity} constant $1$;
\item
the binary \emph{addition} operator ${} + {}$;
\item
the binary \emph{multiplication} operator ${} \mmul {}$;
\item
the unary \emph{additive inverse} operator $- {}$;
\end{itemize}
The signature of inversive meadows consists of the constants and
operators from the signature of commutative rings with a multiplicative
identity and in addition:
\begin{itemize}
\item
the unary \emph{multiplicative inverse} operator ${}\minv$.
\end{itemize}
The signature of divisive meadows consists of the constants and
operators from the signature of commutative rings with a multiplicative
identity and in addition:
\begin{itemize}
\item
the binary \emph{division} operator ${} \mdiv {}$.
\end{itemize}
The signatures of inversive and divisive arithmetical meadows with zero
are the signatures of inversive and divisive meadows with the additive
inverse operator $- {}$ removed.
The signatures of inversive and divisive arithmetical meadows
are the signatures of inversive and divisive arithmetical meadows with
zero with the additive identity constant $0$ removed.
We write:
\begin{ldispl}
\begin{array}{@{}l@{\;}c@{\;}l@{}}
\sigcr    & \mathrm{for} & \set{0,1,{} + {},{} \mmul {}, - {}}\;,
\\
\sigimd   & \mathrm{for} & \sigcr \union \set{{}\minv}\;,
\\
\sigdmd   & \mathrm{for} & \sigcr \union \set{{} \mdiv {}}\;,
\\
\sigiamdz & \mathrm{for} & \sigimd \diff \set{- {}}\;,
\\
\sigdamdz & \mathrm{for} & \sigdmd \diff \set{- {}}\;,
\\
\sigiamd  & \mathrm{for} & \sigiamdz \diff \set{0}\;,
\\
\sigdamd  & \mathrm{for} & \sigdamdz \diff \set{0}\;.
\end{array}
\end{ldispl}

We use infix notation for the binary operators, prefix notation for the
unary operator $- {}$, and postfix notation for the unary operator
${}\minv$.
We use the usual precedence convention to reduce the need for
parentheses.
We denote the numerals $0$, $1$, $1 + 1$, $(1 + 1) + 1$, \ldots~by
$\ul{0}$, $\ul{1}$, $\ul{2}$, $\ul{3}$, \ldots~and we use the notation
$p^n$ for exponentiation with a natural number as exponent.
Formally, we define $\ul{n}$ inductively by $\ul{0} = 0$, $\ul{1} = 1$
and $\ul{n + 2} = \ul{n} + 1$ and we define, for each term $p$ over the
signature of inversive meadows or the signature of divisive meadows,
$p^n$ inductively by $p^0 = 1$ and $p^{n+1} = p^n \mmul p$.

The axioms of a commutative ring with a multiplicative identity are the
equations given in Table~\ref{eqns-commutative-ring}.
We write:
\begin{ldispl}
\begin{array}{@{}l@{\;}c@{\;}l@{}}
\eqnscr  &
\multicolumn{2}{@{}l@{}}
 {\mathrm{for\; the\; set\; of\; all\; equations\; in\; Table\;
          \ref{eqns-commutative-ring}}\;,}
\\
\eqnsacrz & \mathrm{for} &
\eqnscr \diff \set{x + (-x) = 0}\;,
\\
\eqnsacr  & \mathrm{for} &
\eqnsacrz \diff \set{x + 0 = x}\;.
\end{array}
\end{ldispl}
\begin{table}[!t]
\caption{Axioms of a commutative ring with a multiplicative identity}
\label{eqns-commutative-ring}
\begin{eqntbl}
\begin{eqncol}
(x + y) + z = x + (y + z)                                             \\
x + y = y + x                                                         \\
x + 0 = x                                                             \\
x + (-x) = 0
\end{eqncol}
\qquad\quad
\begin{eqncol}
(x \mmul y) \mmul z = x \mmul (y \mmul z)                             \\
x \mmul y = y \mmul x                                                 \\
x \mmul 1 = x                                                         \\
x \mmul (y + z) = x \mmul y + x \mmul z
\end{eqncol}
\end{eqntbl}
\end{table}
The equations in $\eqnsacrz$ are the equations from $\eqnscr$ in which
the additive inverse operator $- {}$ does not occur.
The equations in $\eqnsacr$ are the equations from $\eqnsacrz$ in which
the additive identity constant $0$ does not occur.

To axiomatize inversive and divisive arithmetical meadows, we need
additional equations.
We write:
\begin{ldispl}
\begin{array}{@{}l@{\;}c@{\;}l@{}}
\eqnsiamd & \mathrm{for} &
\eqnsacr \union \set{x \mmul x\minv = 1}\;,
\\
\eqnsdamd & \mathrm{for} &
\eqnsacr \union \set{ x \mdiv x = 1}\;.
\end{array}
\end{ldispl}
The class of \emph{inversive arithmetical meadows} is the class of all
algebras over the signature $\sigiamd$ that satisfy the equations
$\eqnsiamd$ and
the class of \emph{divisive arithmetical meadows} is the class of all
algebras over the signature $\sigdamd$ that satisfy the equations
$\eqnsdamd$.

We state and prove two equational facts about numerals that will be used
in later proofs.
\begin{lemma}
\label{lemma-numerals}
For all $n,m \in \Nat \diff \set{0}$, we have that
$\eqnsacr \deriv \ul{n + m} = \ul{n} + \ul{m}$ and
$\eqnsacr \deriv \ul{n \mmul m} = \ul{n} \mmul \ul{m}$.
\end{lemma}
\begin{proof}
The fact that $\ul{n + m} = \ul{n} + \ul{m}$ is derivable from
$\eqnsacr$ is easily proved by induction on $n$.
The basis step is trivial.
The inductive step goes as follows:
$\ul{(n + 1) + m} = \ul{(n + m) + 1} = \ul{n + m} + 1 =
 \ul{n} + \ul{m} + 1 = \ul{n} + 1 + \ul{m} = \ul{n + 1} + \ul{m}$.
The fact that $\ul{n \mmul m} = \ul{n} \mmul \ul{m}$ is derivable from
$\eqnsacr$ is easily proved by induction on $n$, using that
$\ul{n + m} = \ul{n} + \ul{m}$ is derivable from $\eqnsacr$.
The basis step is trivial.
The inductive step goes as follows:
$\ul{(n + 1) \mmul m} = \ul{n \mmul m + 1 \mmul m} =
 \ul{n \mmul m} + \ul{1 \mmul m} =
 \ul{n} \mmul \ul{m} + \ul{1} \mmul \ul{m} =
 (\ul{n} + \ul{1}) \mmul \ul{m} = \ul{n + 1} \mmul \ul{m}$.
\qed
\end{proof}

We state and prove two useful facts about the multiplicative inverse
operator derivable from $\eqnsiamd$.
\begin{lemma}
\label{lemma-mult-inv-ratia}
We have $\eqnsiamd \deriv (x\minv)\minv = x$ and
$\eqnsiamd \deriv (x \mmul y)\minv = x\minv \mmul y\minv$.
\end{lemma}
\begin{proof}
We derive $(x\minv)\minv = x$ from $\eqnsiamd$ as
follows:
$(x\minv)\minv = 1 \mmul (x\minv)\minv =
 (x \mmul x\minv) \mmul (x\minv)\minv =
 x \mmul (x\minv \mmul (x\minv)\minv) =
 x \mmul 1 = x$.
We derive $(x \mmul y)\minv = x\minv \mmul y\minv$ from $\eqnsiamd$ as
follows:
$(x \mmul y)\minv = 1 \mmul (1 \mmul (x \mmul y)\minv) =\linebreak[2]
 (x \mmul x\minv) \mmul ((y \mmul y\minv) \mmul (x \mmul y)\minv) =
 (x\minv \mmul y\minv) \mmul ((x \mmul y) \mmul (x \mmul y)\minv) =
 (x\minv \mmul y\minv) \mmul 1 = x\minv \mmul y\minv$.
\qed
\end{proof}

To axiomatize inversive and arithmetical meadows with zero, we need
other additional equations.
We write:
\begin{ldispl}
\begin{array}{@{}l@{\;\;}c@{\;\,}l@{}}
\eqnsiamdz & \mathrm{for} &
\eqnsacrz \union \set{{(x\minv)}\minv = x,\;
                      x \mmul (x \mmul x\minv) = x}\;,
\\
\eqnsdamdz & \mathrm{for} &
\eqnsacrz \union \set{1 \mdiv (1 \mdiv x) = x,\;
                      (x \mmul x) \mdiv x = x,\;
                      x \mdiv y = x \mmul (1 \mdiv y)}\;.
\end{array}
\end{ldispl}
The class of \emph{inversive arithmetical meadows with zero} is the
class of all algebras over the signature $\sigiamdz$ that satisfy the
equations $\eqnsiamdz$ and
the class of \emph{divisive arithmetical meadows with zero} is the
class of all algebras over the signature $\sigdamdz$ that satisfy the
equations $\eqnsdamdz$.

We state and prove two useful facts about the additive identity
constant derivable from $\eqnsiamdz$.
\begin{lemma}
\label{lemma-zero-ratiaz}
We have $\eqnsiamdz \deriv 0 \mmul x = 0$ and
$\eqnsiamdz \deriv 0\minv = 0$.
\end{lemma}
\begin{proof}
Firstly, we derive $x + y = x \Implies y = 0$ from $\eqnsiamdz$ as
follows:
$x + y = x \Implies 0 + y = 0 \Implies y + 0 = 0 \Implies y = 0$.
Secondly, we derive $x + 0 \mmul x = x$ from $\eqnsiamdz$ as follows:
$x + 0 \mmul x = x \mmul 1 + 0 \mmul x = 1 \mmul x + 0 \mmul x =
 (1 + 0) \mmul x = 1 \mmul x = x \mmul 1 = x$.
From $x + y = x \Implies y = 0$ and $x + 0 \mmul x = x$, it follows that
$0 \mmul x = 0$.
We derive $0\minv = 0$ from $\eqnsiamdz$ as follows:
$0\minv = 0\minv \mmul (0\minv \mmul (0\minv)\minv) =
 (0\minv)\minv \mmul (0\minv \mmul 0\minv) =
 0 \mmul (0\minv \mmul 0\minv) = 0$.
\qed
\end{proof}

We state and prove a useful fact about the multiplicative inverse
operator derivable from $\eqnsiamdz$.
\begin{lemma}
\label{lemma-mult-inv-ratiaz}
We have $\eqnsiamdz \deriv (x \mmul y)\minv = x\minv \mmul y\minv$.
\end{lemma}
\begin{proof}
Proposition~2.8 from~\cite{BHT09a} states that
$(x \mmul y)\minv = x\minv \mmul y\minv$ is derivable from
$\eqnsiamdz \union \set{x + 0 = x, x + (-x) = 0}$.
The proof of this proposition given in~\cite{BHT09a} goes through
because no use is made of the equations $x + 0 = x$ and $x + (-x) = 0$.
\qed
\end{proof}

We state and prove a fact about the additive identity constant that will
be used in a later proof.
\begin{lemma}
\label{lemma-zero-elim-ratiaz}
For each $\sigiamdz$-term $t$, either $\eqnsiamdz \deriv t = 0$ or there
exists a $\sigiamd$-term $t'$ such that $\eqnsiamdz \deriv t = t'$.
\end{lemma}
\begin{proof}
The proof is easy by induction on the structure of $t$, using
Lemma~\ref{lemma-zero-ratiaz}.
\qed
\end{proof}

\section{Interlude on Meadows}
\label{sect-meadows}

For completeness, we shortly discuss inversive and divisive meadows.

To axiomatize inversive and divisive meadows, we need the following sets
of equations:
\begin{ldispl}
\begin{array}{@{}l@{\;}c@{\;}l@{}}
\eqnsimd & \mathrm{for} &
\eqnscr \union \set{{(x\minv)}\minv = x,\;
                    x \mmul (x \mmul x\minv) = x}\;,
\\
\eqnsdmd & \mathrm{for} &
\eqnscr \union \set{1 \mdiv (1 \mdiv x) = x,\;
                    (x \mmul x) \mdiv x = x,\;
                    x \mdiv y = x \mmul (1 \mdiv y)}\;.
\end{array}
\end{ldispl}
The class of \emph{inversive meadows} is the class of all algebras over
the signature $\sigimd$ that satisfy the equations $\eqnsimd$; and
the class of \emph{divisive meadows} is the class of all algebras over
the signature $\sigdmd$ that satisfy the equations $\eqnsdmd$.

A meadow is \emph{non-trivial} if it satisfies $0 \neq 1$; and
a meadow is a \emph{cancellation meadow} if it satisfies
$x \neq 0 \And x \mmul y = x \mmul z \Implies y = z$.
In~\cite{BT07a}, an inversive cancellation meadow is called a
\emph{zero-totalized field}.

Recently, we found out in~\cite{Ono83a} that inversive meadows were
already introduced by Komori~\cite{Kom75a} in a report from 1975, where
they go by the name of \emph{desirable pseudo-fields}.
In~\cite{Ono83a}, we also found an axiomatization of inversive meadows
which differs from the one given above.
We came across this paper by a reference in~\cite{Pie92a}.
We propose the names \emph{inversive Komori field} and \emph{divisive
Komori field} as alternatives for inversive non-trivial cancellation
meadow and divisive non-trivial cancellation meadow, respectively.

$\Ratzi$, the inversive meadow of rational numbers, is defined as
follows:
\begin{ldispl}
\Ratzi =
I(\sigimd,\eqnsimd \union
          \set{(1 + x^2 + y^2) \mmul (1 + x^2 + y^2)\minv = 1})\;.
\end{ldispl}
$\Ratzd$, the divisive meadow of rational numbers, is defined as
follows:
\begin{ldispl}
\Ratzd =
I(\sigdmd,\eqnsdmd \union
          \set{(1 + x^2 + y^2) \mdiv (1 + x^2 + y^2) = 1})\;.
\end{ldispl}
$\Ratzi$ is an inversive Komori field and $\Ratzd$ is a divisive Komori
field.

Moss showed in~\cite{Mos01a} that there exists an initial algebra
specification of $\Rat$ with just one hidden function.
The initial algebra specification of $\Ratzi$ given above is without
hidden functions.

Inversive meadows have been extended with signum, floor and ceiling
operations in~\cite{BP08a}, differentiation operations in~\cite{BP08b},
and a square root operation in~\cite{BB09a}.

\section{Arithmetical Meadows of Rational Numbers}
\label{sect-arith-meadows-rat}

We obtain inverse and divisive arithmetical meadows closely related to
the field of rational numbers as the initial algebras of equational
specifications.
As usual, we write $I(\Sigma,E)$ for the initial algebra among the
algebras over the signature $\Sigma$ that satisfy the equations $E$
(see e.g.~\cite{BT87a}).

$\Ratia$, the inversive arithmetical meadow of rational numbers, is
defined as follows:
\begin{ldispl}
\Ratia = I(\sigiamd,\eqnsiamd)\;.
\end{ldispl}
$\Ratda$, the divisive arithmetical meadow of rational numbers, is
defined as follows:
\begin{ldispl}
\Ratda = I(\sigdamd,\eqnsdamd)\;.
\end{ldispl}
$\Ratia$ and $\Ratda$ are the initial algebras in the class of inversive
arithmetical meadows and the class of divisive arithmetical meadows,
respectively.

$\Ratia$ is a subalgebra of a reduct of $\Ratzi$.
\begin{theorem}
\label{theorem-init-alg-spec-ratia}
$\Ratia$ is the subalgebra of the $\sigiamd$-reduct of $\Ratzi$ whose
domain is the set of all positive rational numbers.
\end{theorem}
\begin{proof}
Like in the case of Theorem~3.1 from~\cite{BT07a}, it is sufficient to
prove that, for each closed term $t$ over the signature $\sigiamd$,
there exists a unique term $t'$ in the set
\begin{ldispl}
\set{\ul{n} \mmul \ul{m}\minv \where
n,m \in \Nat \diff \set{0} \And \nm{gcd}(n,m) = 1}
\end{ldispl}
such that $\eqnsiamd \deriv t = t'$.
Like in the case of Theorem~3.1 from~\cite{BT07a}, this is proved by
induction on the structure of $t$, using Lemmas~\ref{lemma-numerals}
and~\ref{lemma-mult-inv-ratia}.
The proof is similar, but simpler owing to:
(i)~the absence of terms of the forms $0$ and $-t'$;\linebreak[2]
(ii)~the absence of terms of the forms $\ul{0}$ and
$- (\ul{n} \mmul \ul{m}\minv)$ among the terms that exist by the
induction hypothesis;
(iii)~the presence of the axiom $x \mmul x\minv = 1$.
\qed
\end{proof}
The fact that $\Ratda$ is the initial algebra in the class of
divisive arithmetical meadows is proved similarly.

Derivability of equations from the equations of the initial algebra
specifications of $\Ratia$ and $\Ratda$ is decidable.
\begin{theorem}
\label{theorem-decidability-ratia}
For all $\sigiamd$-terms $t$ and $t'$, it is decidable whether
$\eqnsiamd \deriv t = t'$.
\end{theorem}
\begin{proof}
For each $\sigiamd$-term $r$, there exist $\sigiamd$-terms $r_1$ and
$r_2$ in which the multiplicative inverse operator do not occur such
that $\eqnsiamd \deriv r = r_1 \mmul r_2\minv$.
The proof of this fact is easy by induction on the structure of $r$,
using Lemma~\ref{lemma-mult-inv-ratia}.
Inspection of the proof yields that there is an effective way to find
witnessing terms.

For each closed $\sigiamd$-term $r$ in which the multiplicative inverse
operator does not occur there exists a $k \in \Nat \diff \set{0}$, such
that $\eqnsiamd \deriv r = \ul{k}$.
The proof of this fact is easy by induction on the structure of $r$.
Moreover, for each $\sigiamd$-term $r$ in which the multiplicative
inverse operator does not occur there exists a $\sigiamd$-term $r'$ of
the form
$\sum_{i_1=1}^{n_1} \ldots \sum_{i_m=1}^{n_m}
  \ul{k_{i_1 \ldots i_m}} \mmul
  x_1^{i_1} \mmul {} \cdots {} \mmul x_m^{i_m}$,
where $k_{i_1 \ldots i_m} \in \Nat \diff \set{0}$ for each
$i_1 \in [1,n_1]$, \ldots, $i_m \in [1,n_m]$ and $x_1,\ldots,x_m$ are
variables, such that $\eqnsiamd \deriv r = r'$.
The proof of this fact is easy by induction on the structure of $r$,
using the previous fact.
Inspection of the proof yields that there is an effective way to find
a witnessing term.
Terms of the form described above are polynomials in several variables
with positive integer coefficients.

Let $t_1$, $t_2$, $t'_1$, $t'_2$ be $\sigiamd$-terms in which the
multiplicative inverse operator do not occur such that
$\eqnsiamd \deriv t = t_1 \mmul {t_2}\minv$ and
$\eqnsiamd \deriv t' = t'_1 \mmul {t'_2}\minv$.
Moreover, let $s$ and $s'$ be $\sigiamd$-terms of the form
$\sum_{i_1=1}^{n_1} \ldots \sum_{i_m=1}^{n_m}
  \ul{k_{i_1 \ldots i_m}} \mmul
  x_1^{i_1} \mmul {} \cdots {} \mmul x_m^{i_m}$,
where $k_{i_1 \ldots i_m} \in \Nat \diff \set{0}$ for each
$i_1 \in [1,n_1]$, \ldots, $i_m \in [1,n_m]$ and $x_1,\ldots,x_m$ are
variables, such that $\eqnsiamd \deriv t_1 \mmul t'_2 = s$ and
$\eqnsiamd \deriv t'_1 \mmul t_2 = s'$.
We have that
$\eqnsiamd \deriv t = t'$ iff
$\eqnsiamd \deriv t_1 \mmul {t_2}\minv = t'_1 \mmul {t'_2}\minv$ iff
$\eqnsiamd \deriv t_1 \mmul t'_2 = t'_1 \mmul t_2$ iff
$\eqnsiamd \deriv s = s'$.
Moreover, we have that $\eqnsiamd \deriv s = s'$ only if $s$ and $s'$
denote the same function on positive real numbers in the inversive
arithmetical meadow of positive real numbers.
The latter is decidable because polynomials in several variables with
positive integer coefficients denote the same function on positive real
numbers in the inversive arithmetical meadow of positive real numbers
only if they are syntactically equal.
\qed
\end{proof}
The fact that derivability of equations from the equations of the
initial algebra specification of $\Ratda$ is decidable is proved
similarly.

We obtain inverse and divisive arithmetical meadows with zero closely
related to the field of rational numbers as the initial algebras of
equational specifications.

$\Ratiaz$, the inversive arithmetical meadow of rational numbers with
zero, is defined as follows:
\begin{ldispl}
\Ratiaz =
I(\sigiamdz,\eqnsiamdz \union
            \set{(1 + x^2 + y^2) \mmul (1 + x^2 + y^2)\minv = 1})\;.
\end{ldispl}
$\Ratdaz$, the divisive arithmetical meadow of rational numbers with
zero, is defined as follows:
\begin{ldispl}
\Ratdaz =
I(\sigdamdz,\eqnsdamdz \union
            \set{(1 + x^2 + y^2) \mdiv (1 + x^2 + y^2) = 1})\;.
\end{ldispl}

$\Ratiaz$ is a subalgebra of a reduct of $\Ratzi$.
First we prove a fact that is useful in the proving this result.
\begin{lemma}
\label{lemma-mult-inv-exist-ratiaz}
It follows from
$\eqnsiamdz \union \set{(1 + x^2 + y^2) \mmul (1 + x^2 + y^2)\minv = 1}$
that $\ul{n}$ has a multiplicative inverse for each
$n \in \Nat \diff \set{0}$.
\end{lemma}
\begin{proof}
In the proof of Theorem~3 from~\cite{BM09g}, it is among other things
proved that it follows from
$\eqnsiamdz \union \set{x + (-x) = 0} \union
   \set{(1 + x^2 + y^2) \mmul (1 + x^2 + y^2)\minv = 1}$
that $\ul{n}$ has a multiplicative inverse for each
$n \in \Nat \diff \set{0}$.
The proof concerned goes through because no use is made of the equation
$x + (-x) = 0$.
\qed
\end{proof}
\begin{theorem}
\label{theorem-init-alg-spec-ratiaz}
$\Ratiaz$ is the subalgebra of the $\sigiamdz$-reduct of $\Ratzi$ whose
domain is the set of all non-negative rational numbers.
\end{theorem}
\begin{proof}
Like in the case of Theorem~\ref{theorem-init-alg-spec-ratia}, it is
sufficient to prove that, for each closed term $t$ over the signature
$\sigiamd$, there exists a unique term $t'$ in the set
\begin{ldispl}
\set{\ul{0}} \union
\set{\ul{n} \mmul \ul{m}\minv \where
     n,m \in \Nat \diff \set{0} \And \nm{gcd}(n,m) = 1}
\end{ldispl}
such that
$\eqnsiamdz \union
 \set{(1 + x^2 + y^2) \mmul (1 + x^2 + y^2)\minv = 1} \deriv t = t'$.
Like in the case of Theorem~\ref{theorem-init-alg-spec-ratia}, this is
proved by induction on the structure of $t$, now using
Lemmas~\ref{lemma-numerals}, \ref{lemma-zero-ratiaz}
and~\ref{lemma-mult-inv-ratiaz}.
The proof is similar, but more complicated owing to:
(i)~the presence of terms of the form $0$;
(ii)~the presence of terms of the form $\ul{0}$ among the terms that
exist by the induction hypothesis;
(iii)~the absence of the axiom $x \mmul x\minv = 1$.
Because of the last point, use is made of
Lemma~\ref{lemma-mult-inv-exist-ratiaz}.
\qed
\end{proof}
The fact that $\Ratdaz$ is the initial algebra in the class of divisive
arithmetical meadows with zero is proved similarly.
\pagebreak[2]

An alternative initial algebra specification of $\Ratiaz$ is obtained if
the equation $(1 + x^2 + y^2) \mmul (1 + x^2 + y^2)\minv = 1$ is
replaced by
$(x \mmul (x + y)) \mmul (x \mmul (x + y))\minv = x \mmul x\minv$.
\begin{theorem}
\label{theorem-alt-init-alg-spec-ratiaz}
$\Ratiaz \cong
 I(\sigiamdz,
   \eqnsiamdz \union
   \set{(x \mmul (x + y)) \mmul (x \mmul (x + y))\minv =
        x \mmul x\minv})$.
\end{theorem}
\begin{proof}
It is sufficient to prove that
$(x \mmul (x + y)) \mmul (x \mmul (x + y))\minv = x \mmul x\minv$
is valid in $\Ratiaz$ and
$(1 + x^2 + y^2) \mmul (1 + x^2 + y^2)\minv = 1$ is valid in
$I(\sigiamdz,
   \eqnsiamdz \union
   \set{(x \mmul (x + y)) \mmul (x \mmul (x + y))\minv =
        x \mmul x\minv})$.
It follows from Lemma~\ref{lemma-mult-inv-ratiaz}, and the
associativity and commutativity of ${}\mmul{}$, that
$(x \mmul (x + y)) \mmul (x \mmul (x + y))\minv = x \mmul x\minv \Iff
\linebreak[2]
 (x \mmul x\minv) \mmul ((x + y) \mmul (x + y)\minv) =
 x \mmul x\minv$ is derivable from $\eqnsiamdz$.
This implies that
$(x \mmul (x + y)) \mmul (x \mmul (x + y))\minv = x \mmul x\minv$
is valid in $\Ratiaz$ iff
$(x \mmul x\minv) \mmul ((x + y) \mmul (x + y)\minv) = x \mmul x\minv$
is valid in $\Ratiaz$.
The latter is easily established by distinction
between the cases $x = 0$ and $x \neq 0$.
To show that
$(1 + x^2 + y^2) \mmul (1 + x^2 + y^2)\minv = 1$ is valid in
$I(\sigiamdz,
   \eqnsiamdz \union
   \set{(x \mmul (x + y)) \mmul (x \mmul (x + y))\minv =
        x \mmul x\minv})$,
it is sufficient to derive
$(1 + x^2 + y^2) \mmul (1 + x^2 + y^2)\minv = 1$ from
$\eqnsiamdz \union
 \set{(x \mmul x\minv) \mmul ((x + y) \mmul (x + y)\minv) =
      x \mmul x\minv}$.
The derivation is fully trivial with the exception of the first step,
viz.\ substituting $1$ for $x$ and $x^2 + y^2$ for $y$ in
$(x \mmul x\minv) \mmul ((x + y) \mmul (x + y)\minv) = x \mmul x\minv$.
\qed
\end{proof}
An alternative initial algebra specification of $\Ratdaz$ is obtained in
the same vein.

In $\Ratiaz$, the \emph{general inverse law}
$x \neq 0 \Implies x \mmul x\minv = 1$ is valid.
Derivability of equations from the equations of the alternative initial
algebra specification of $\Ratiaz$ and the general inverse law is
decidable.
First we prove a fact that is useful in proving this decidability
result.
\begin{lemma}
\label{lemma-derivability-gil}
For all $\sigiamd$-terms $t$ in which no other variables than
$x_1,\ldots,x_n$ occur, 
$\eqnsiamdz \union
 \set{(x \mmul (x + y)) \mmul (x \mmul (x + y))\minv =
      x \mmul x\minv} \union
 \set{x_1 \mmul x_1\minv = 1,\ldots,\linebreak[2]x_n \mmul x_n\minv = 1}
  \deriv_{x_1,\ldots,x_n} t \mmul t\minv = 1$.
\end{lemma}
\begin{proof}
The proof is easy by induction on the structure of $t$, using
Lemma~\ref{lemma-mult-inv-ratiaz}.
\qed
\end{proof}

\begin{theorem}
\label{theorem-decidability-ratiaz}
For all $\sigiamdz$-terms $t$ and $t'$, it is decidable whether
$\eqnsiamdz \union \linebreak
 \set{(x \mmul (x + y)) \mmul (x \mmul (x + y))\minv =
      x \mmul x\minv} \union
 \set{x \neq 0 \Implies x \mmul x\minv = 1} \deriv t = t'$.
\end{theorem}
\begin{proof}
Let $\eqnsiamdzx =
 \eqnsiamdz \union
 \set{(x \mmul (x + y)) \mmul (x \mmul (x + y))\minv =
      x \mmul x\minv} \union \linebreak
 \set{x \neq 0 \Implies x \mmul x\minv = 1}$.
We prove that $\eqnsiamdzx \deriv t = t'$ is decidable by induction on
the number of variables occurring in $t = t'$.
In the case where the number of variables is $0$, we have that
$\eqnsiamdzx \deriv t = t'$ iff $\Ratiaz \models t = t'$ iff
$\eqnsiamdz \union
 \set{(1 + x^2 + y^2) \mmul (1 + x^2 + y^2)\minv = 1} \deriv t = t'$.
The last is decidable because, by the proof of
Theorem~\ref{theorem-init-alg-spec-ratiaz}, there exist
unique terms $s$ and $s'$ in the set
$\set{\ul{0}} \union
 \set{\ul{n} \mmul \ul{m}\minv \where
      n,m \in \Nat \diff \set{0} \And \nm{gcd}(n,m) = 1}$
such that
$\eqnsiamdz \union
 \set{(1 + x^2 + y^2) \mmul (1 + x^2 + y^2)\minv = 1} \deriv t = s$ and
$\eqnsiamdz \union
 \set{(1 + x^2 + y^2) \mmul \linebreak[2] (1 + x^2 + y^2)\minv = 1}
  \deriv t' = s'$,
and inspection of that proof yields that there is an effective way to
find $s$ and $s'$.
Hence, in the case where the number of variables is $0$,
$\eqnsiamdzx \deriv t = t'$ is decidable.
In the case where the number of variables is $n + 1$, suppose that the
variables are $x_1,\ldots,x_{n+1}$.
Let $s$ be such that $\eqnsiamdz \deriv t = s$ and $s$ is either a
$\sigiamd$-term or the constant $0$ and let $s'$ be such that
$\eqnsiamdz \deriv t' = s'$ and $s'$ is either a $\sigiamd$-term or the
constant $0$.
Such $s$ and $s'$ exist by Lemma~\ref{lemma-zero-elim-ratiaz}, and
inspection of the proof of that lemma yields that there is an effective
way to find $s$ and $s'$.
We have that $\eqnsiamdzx \deriv t = t'$ iff
$\eqnsiamdzx \deriv s = s'$.
In the case where not both $s$ and $s'$ are $\sigiamd$-terms,
$\eqnsiamdzx \deriv s = s'$ only if $s$ and $s'$ are syntactically
equal.
Hence, in this case, $\eqnsiamdzx \deriv t = t'$ is decidable.
In the case where both $s$ and $s'$ are $\sigiamd$-terms, by the general
inverse law, we have that
$\eqnsiamdzx \deriv s = s'$ iff
$\eqnsiamdzx \deriv s\subst{0}{x_i} =
 s'\subst{0}{x_i}$
for all $i \in [1,n + 1]$ and
$\eqnsiamdzx \union
 \set{x_1 \mmul x_1\minv = 1,\ldots,x_{n+1} \mmul x_{n+1}\minv = 1}
  \deriv_{x_1,\ldots,x_{n+1}} s = s'$.
By Lemma~\ref{lemma-derivability-gil}, we have that
$\eqnsiamdzx \union
 \set{x_1 \mmul x_1\minv = 1,\ldots,x_{n+1} \mmul x_{n+1}\minv = 1}
  \deriv_{x_1,\ldots,x_{n+1}} s = s'$ iff
$\eqnsiamd \deriv s = s'$.
For each $i \in [1,n + 1]$,
$\eqnsiamdzx \deriv s\subst{0}{x_i} =
 s'\subst{0}{x_i}$
is decidable because the number of variables occurring in
$s\subst{0}{x_i} = s'\subst{0}{x_i}$ is $n$.
Moreover, we know from Theorem~\ref{theorem-decidability-ratia} that
$\eqnsiamd \deriv s = s'$ is decidable.
Hence, in the case where both $s$ and $s'$ are $\sigiamd$-terms,
$\eqnsiamdzx \deriv t = t'$ is decidable as well.
\qed
\end{proof}
The fact that derivability of equations from the equations of the
alternative initial algebra specification of $\Ratdaz$ and
$x \neq 0 \Implies x \mdiv x = 1$ is decidable is proved similarly.
We remark that it is an open problem whether derivability of equations
from the equations of the alternative initial algebra specifications of
$\Ratiaz$ and $\Ratdaz$ is decidable.

\section{Outstanding Questions about Arithmetical Meadows}
\label{sect-questions}

The following are some outstanding questions with regard to arithmetical
meadows:
\begin{enumerate}
\item
Is the initial algebra specification of $\Ratzi$ a conservative
extension of the initial algebra specifications of $\Ratia$ and
$\Ratiaz$?
\item
Do $\Ratia$ and $\Ratiaz$ have initial algebra specifications that
constitute complete term rewriting systems (modulo associativity and
commutativity of ${}+{}$ and~${}\mmul{}$)?
\item
Do $\Ratia$ and $\Ratiaz$ have $\omega$-complete initial algebra
specifications?
\item
What are the complexities of derivability of equations from
$\eqnsiamd$ and
$\eqnsiamdz \union
 \set{(x \mmul (x + y)) \mmul (x \mmul (x + y))\minv =
      x \mmul x\minv,\;
      x \neq 0 \Implies x \mmul x\minv = 1}$?
\item
Is derivability of equations from
$\eqnsiamdz \union
 \set{(x \mmul (x + y)) \mmul (x \mmul (x + y))\minv =
      x \mmul x\minv} \deriv t = t'$ decidable?
\item
Do we have
$\Ratiaz \cong
 I(\sigiamdz,
   \eqnsiamdz \union \set{(1 + x^2) \mmul (1 + x^2)\minv = 1})$?
\end{enumerate}
These questions are formulated for the inversive case, but they have
counterparts for the divisive case of which some might lead to different
answers.

\section{Partial Arithmetical Meadows with Zero}
\label{sect-partial-arith-meadows}

Following~\cite{BM09g}, we introduce in this section simple
constructions of partial inversive and divisive arithmetical meadows
with zero from total ones.

We take the position that partial algebras should be made from total
ones.
For the case that we are engaged in, this means that relevant partial
arithmetical meadows with zero are obtained from arithmetical meadows
with zero by making certain operations undefined for certain arguments.

Let $\Mdiaz$ be an inversive arithmetical meadow with zero.
Then it makes sense to construct one partial inversive arithmetical
meadow with zero from $\Mdiaz$:
\begin{itemize}
\item
$0\minv \undef \Mdiaz$ is the partial algebra that is obtained from
$\Mdiaz$ by making $0\minv$ undefined.
\end{itemize}
Let $\Mddaz$ be a divisive arithmetical meadow with zero.
Then it makes sense to construct two partial divisive arithmetical
meadows with zero from $\Mddaz$:
\begin{itemize}
\item
$\Quant \mdiv 0 \undef \Mddaz$ is the partial algebra that is obtained
from $\Mddaz$ by making\linebreak[2] $q \mdiv 0$ undefined for all $q$
in the domain of $\Mddaz$;
\item
$(\Quant \diff \set{0}) \mdiv 0 \undef \Mddaz$ is the partial algebra
that is obtained from $\Mddaz$ by\linebreak[2] making $q \mdiv 0$
undefined for all $q$ in the domain of $\Mddaz$ different from $0$.
\end{itemize}

Clearly, the partial arithmetical meadow constructions are special cases
of a more general partial algebra construction for which we have coined
the term \emph{punching}.
Presenting the details of the general construction is outside the scope
of the current paper.

The partial arithmetical meadow constructions described above yield the
\linebreak[2] following three partial arithmetical meadows with zero
related to rational
numbers:
\begin{ldispl}
\begin{eqncol}
0\minv \undef \Ratiaz\;,
\qquad
\Quant \mdiv 0 \undef \Ratdaz\;,
\qquad
(\Quant \diff \set{0}) \mdiv 0 \undef \Ratdaz\;.
\end{eqncol}
\end{ldispl}
These algebras have been obtained by means of the well-known initial
algebra construction and a simple partial algebra construction.
The merits of this approach are discussed in~\cite{BM09g}.

At first sight, the absence of the additive inverse operator does not
seem to add anything new to the treatment of punched meadows
in~\cite{BM09g}.
However, this is not quite the case.
Consider $0\minv \undef \Ratiaz$.
In the case of this algebra, there is a useful syntactic criterion for
``being defined''.
The set $\Def$ of defined terms and the auxiliary set $\Nz$ of non-zero
terms can be inductively defined by:
\begin{itemize}
\item
$1 \in \Nz$;
\item
if $x \in \Nz$, then $x + y \in \Nz$ and $y + x \in \Nz$;
\item
if $x \in \Nz$ and $y \in \Nz$, then $x \mmul y \in \Nz$;
\item
if $x \in \Nz$, then $x\minv \in \Nz$;
\item
$0 \in \Def$;
\item
if $x \in \Nz$, then $x \in \Def$;
\item
if $x \in \Def$ and $y \in \Def$, then $x + y \in \Def$ and
$x \mmul y \in \Def$.
\end{itemize}
This indicates that the absence of the additive inverse operator allows
a typing based solution to problems related to ``division by zero'' in
elementary school mathematics.
So there may be a point in dealing first and thoroughly with
non-negative rational numbers in a setting where division by zero is not
defined.

Working in $\Ratia$ simplifies matters even more because there is no
distinction between terms and defined terms.
Again, this may be of use in the teaching of mathematics at elementary
school.

\section{Arithmetical Meadows and Regular Arithmetical Rings}
\label{sect-relation-arith-rings}

We can define commutative arithmetical rings with a multiplicative
identity in the same vein as arithmetical meadows.
Moreover, we can define commutative von Neumann regular arithmetical
rings with a multiplicative identity as commutative arithmetical rings
with a multiplicative identity satisfying the regularity condition
$\Forall{x}{\Exists{y}{x \mmul (x \mmul y) = x}}$.

The following theorem states that commutative von Neumann regular
arithmetical rings with a multiplicative identity are related to
inversive arithmetical meadows like commutative von Neumann regular
rings with a multiplicative identity are related to inversive meadows.
\begin{theorem}
Each commutative von Neumann regular arithmetical ring with a
multiplicative identity can be expanded to an inversive arithmetical
meadow, and this expansion is unique.
\end{theorem}
\begin{proof}
Lemma~2.11 from~\cite{BHT09a} states that each commutative von Neumann
regular ring with a multiplicative identity element can be expanded to
an inversive meadow, and this expansion is unique.
The only use that is made of the equations $x + 0 = x$ and
$x + (-x) = 0$ in the proof of this lemma given in~\cite{BHT09a}
originates from the proof of another lemma that is used in the proof.
However, the latter lemma, Lemma~2.12 from~\cite{BHT09a}, concerns the
same property as Proposition~2.3 from~\cite{BRS09a} and in the proof of
this proposition given in~\cite{BRS09a} no use is made of the equations
$x + 0 = x$ and $x + (-x) = 0$.
Hence, there is an alternative proof of Lemma~2.11 from~\cite{BHT09a}
that goes through for the arithmetical case.
\qed
\end{proof}

We can also define commutative arithmetical rings with additive and
multiplicative identities and commutative von Neumann regular
arithmetical rings with additive and multiplicative identities in the
obvious way.
We also have that commutative von Neumann regular arithmetical rings
with additive and multiplicative identities are related to inversive
arithmetical meadows with zero like commutative von Neumann regular
rings with a multiplicative identity are related to inversive meadows.

\section{Conclusions}
\label{sect-conclusions}

We have complemented the signatures of inversive and divisive meadows
with arithmetical versions, and provided equational axiomatizations of
several classes of arithmetical meadows and instances of them related to
the field of rational numbers.
We have answered a number of questions about these classes and instances
of arithmetical meadows, and stated a number of outstanding questions
about them.
In addition, we have discussed partial variants of the instances in
question and an arithmetical version of a well-known mathematical
structure closely related to inversive meadows, namely von Neumann
regular rings.

We remark that the name arithmetical algebra is not always used in the
same way as Peacock~\cite{Pea1830a} used it.
It is sometimes difficult to establish whether the notion in question is
related to Peacock's notion of arithmetical algebra.
For example, it is not clear to us whether the notion of arithmetical
algebra defined in~\cite{Pix72a} is related to Peacock's notion of
arithmetical algebra.

\bibliographystyle{splncs03}
\bibliography{MD}

\begin{thebibliography}{10}
\providecommand{\url}[1]{\texttt{#1}}
\providecommand{\urlprefix}{URL }

\bibitem{BB09a}
Bergstra, J.A., Bethke, I.: Square root meadows. {\tt arXiv:0901.4664v1
  [cs.LO]} at {\tt http://arxiv.org/} (January 2009)

\bibitem{BHT09a}
Bergstra, J.A., Hirshfeld, Y., Tucker, J.V.: Meadows and the equational
  specification of division. Theoretical Computer Science  410(12--13),
  1261--1271 (2009)

\bibitem{BM09g}
Bergstra, J.A., Middelburg, C.A.: Inversive meadows and divisive meadows.
  {\tt arXiv:0907.0540v2 [math.RA]} at {\tt http://arxiv.org/} (July 2009) 

\bibitem{BP08b}
Bergstra, J.A., Ponse, A.: Differential meadows. {\tt arXiv:0804.3336v1
  [math.RA]} at {\tt http://arxiv.org/} (April 2008)

\bibitem{BP08a}
Bergstra, J.A., Ponse, A.: A generic basis theorem for cancellation meadows.
  {\tt arXiv:0803.3969v2 [math.RA]} at {\tt http://arxiv.org/} (March 2008)

\bibitem{BT87a}
Bergstra, J.A., Tucker, J.V.: Algebraic specifications of computable and
  semicomputable data types. Theoretical Computer Science  50(2),  137--181
  (1987)

\bibitem{BT95a}
Bergstra, J.A., Tucker, J.V.: Equational specifications, complete term
  rewriting, and computable and semicomputable algebras. Journal of the ACM
  42(6),  1194--1230 (1995)

\bibitem{BT07a}
Bergstra, J.A., Tucker, J.V.: The rational numbers as an abstract data type.
  Journal of the ACM  54(2),  Article 7 (2007)

\bibitem{BRS09a}
Bethke, I., Rodenburg, P.H., Sevenster, A.: The structure of finite meadows.
  {\tt arXiv:0903.1196v1 [cs.LO]} at {\tt http://arxiv.org/} (March 2009)

\bibitem{Goo79a}
Goodearl, K.R.: {Von Neumann} Regular Rings. Pitman, London (1979)

\bibitem{Kle98a}
Kleiner, I.: A historically focused course in abstract algebra. Mathematics
  Magazine  71(2),  105--111 (1998)

\bibitem{Kom75a}
Komori, Y.: Free algebras over all fields and pseudo-fields. Report 10,
  pp.~9--15, Faculty of Science, Shizuoka University (1975)

\bibitem{McC64a}
McCoy, N.H.: The Theory of Rings. Macmillan, London (1964)

\bibitem{Mos01a}
Moss, L.: Simple equational specifications of rational arithmetics. Discrete
  Mathematics and Theoretical Computer Science  4,  291--300 (2001)

\bibitem{Ono83a}
Ono, H.: Equational theories and universal theories of fields. Journal of the
  Mathematical Society of Japan  35(2),  289--306 (1983)

\bibitem{Pea1830a}
Peacock, G.: A Treatise on Algebra. J. \& J.J. Deighton, Cambridge (1830)

\bibitem{Pie92a}
Pierce, R.S.: Minimal regular rings. In: Contemporary Mathematics, vol. 130,
  pp. 335--348. American Mathematical Society (1992)

\bibitem{Pix72a}
Pixley, A.F.: Completeness in arithmetical algebras. Algebra Universalis  2(1),
   179--196 (1972)

\end{thebibliography}


\end{document}